\documentclass[fleqn]{mat01}
\usepackage{times,mathtimy,amssymb,latexsym,epsfig}
\begin{document}

\setcounter{page}{161}
\firstpage{161}

\makeatletter
\def\artpath#1{\def\@artpath{#1}}
\makeatother \artpath{C:/mathsci-arxiv/may2006/texfiles}

\def\d{\mbox{\rm d}}
\def\e{\mbox{\rm e}}
\def\Rep{\mbox{\rm Rep}}
\def\srad{\mbox{\rm srad}}

\def\mand{\quad \mbox{and} \quad}
\def\cN{{\mathcal N}}
\def\cU{{\mathcal U}}
\def\cF{{\mathcal F}}

\newcommand{\R}{{\mathbb R}}
\newcommand{\N}{{\mathbb N}}
\newcommand{\Z}{{\mathbb Z}}
\newcommand{\lL}{{\mathbb L}}

\newtheorem{theora}{Theorem}
\renewcommand\thetheora{\arabic{section}.\arabic{theora}}
\newtheorem{lem}[theora]{Lemma}
\newtheorem{coro}[theora]{\rm COROLLARY}

\newtheorem{theore}{\bf Theorem}
\newtheorem{case}{\it Case}
\def\notat{\trivlist \item[\hskip \labelsep{{\it Notation}.}]}
\def\claim{\trivlist \item[\hskip \labelsep{\it Claim.}]}

\newtheorem{theoree}{Theorem}
\renewcommand\thetheoree{\it\Roman{theoree}}
\newtheorem{step}[theoree]{Step}

\title{Enveloping $\pmb{\sigma}$-$\pmb{C^{*\!}}$-algebra of a smooth
Frechet algebra\\ crossed product by $\pmb{\R}$, $\pmb{K}$-theory
and differential structure\\ in \hbox{$\pmb{C^{*\!}}$-algebras}}

\markboth{Subhash J Bhatt}{Enveloping $\sigma$-$C^{*\!}$-algebra of a
smooth Frechet algebra}

\author{SUBHASH J BHATT}

\address{Department of Mathematics, Sardar Patel University, Vallabh
Vidyanagar~388~120, India\\
\noindent E-mail: subhashbhaib@yahoo.co.in}

\volume{116}

\mon{May}

\parts{2}

\pubyear{2006}

\Date{MS received 9 August 2004}

\begin{abstract}
Given an $m$-tempered strongly continuous action $\alpha$ of $\R$
by continuous $^{*}$-automorphisms of a Frechet $^{*}$-algebra
$A$, it is shown that the enveloping
\hbox{$\sigma$-$C^{*\!}$-algebra} $E(S(\R,A^{\infty},\alpha))$ of
the smooth Schwartz crossed product $S(\R,A^{\infty},\alpha)$ of
the Frechet algebra $A^{\infty}$ of $C^{\infty}$-elements of $A$
is isomorphic to the \hbox{$\sigma$-$C^{*\!}$-crossed} product
$C^{*}(\R,E(A),\alpha)$ of the enveloping
$\sigma$-$C^{*\!}$-algebra $E(A)$ of $A$ by the induced action.
When $A$ is a hermitian $Q$-algebra, one gets $K$-theory
isomorphism $RK_{*}(S(\R,A^{\infty},\alpha)) =
K_{*}(C^{*}(\R,E(A),\alpha)$ for the representable $K$-theory of
Frechet algebras. An application to the differential structure of
a $C^{*\!}$-algebra defined by densely defined differential
seminorms is given.
\end{abstract}

\keyword{Frechet $^{*}$-algebra; enveloping $\sigma$-$C^{*\!}$-algebra;
smooth crossed product; $m$-tempered action; $K$-theory; differential
structure in $C^{*\!}$-algebras.}

\maketitle

\section{Introduction}

Given a strongly continuous action $\alpha$ of $\R$ by continuous
$^{*}$-automorphisms of a Frechet $^{*}$-algebra $A$, several
crossed product Frechet algebras can be constructed \cite{11,14}.
They include the smooth Schwartz crossed product $S(\R,A,\alpha)$,
the $L^{1}$-crossed products $L^{1}(\R,A,\alpha)$ and
$L^{1}_{|\cdot|}(\R,A,\alpha)$, and the $\sigma$-$C^{*\!}$-crossed
product $C^{*}(\R,A,\alpha)$. Let $E(A)$ denote the enveloping
$\sigma$-$C^{*\!}$-algebra of $A$ \cite{1,6}; and
$(A^{\infty},\tau)$ denote the Frechet $^{*}$-algebra consisting
of all $C^{\infty}$-elements of $A$ with the $C^{\infty}$-topology
$\tau$ (\cite{14}, Appendix~I). The following theorem shows that
for a smooth action, the eveloping algebra of smooth crossed
product is the continuous crossed product of the enveloping
algebra.

\begin{theore}[\!]
Let $\alpha$ be an $m$-tempered strongly continuous action of $\R$
by continuous \hbox{$^{*}$-automorphisms} of a Frechet
$^{*}$-algebra $A$. Let $A$ admit a bounded approximate identity
which is contained in $A^{\infty}$ and which is a bounded
approximate identity for the Frechet algebra $A^{\infty}$. Then
$E(S(\R,A^{\infty},\alpha)) \cong
E(L^{1}_{|\cdot|}(\R,A^{\infty},\alpha))\cong C^{*}(\R,
E(A),\alpha)$. Further{\rm ,} if $\alpha$ is isometric{\rm ,} then
$E(L^{1}(\R,A,\alpha))\cong C^{*}(\R, E(A),\alpha)$.
\end{theore}

Notice that neither $L^{1}(\R,A,\alpha)$ nor
$S(\R,A^{\infty},\alpha)$ need be a subalgebra of
$C^{*}(\R,E(A),\break \alpha)$. A particular case of Theorem~1
when $A$ is a dense subalgebra of $C^{*\!}$-algebra has been
treated in \cite{2}. Let $RK_{*}$ (respectively $K_{*}$) denote
the representable $K$-theory functor (respectively $K$-theory
functor) on Frechet algebras \cite{10}. We have the following
isomorphism of $K$-theory, obtained without direct appeal to
spectral invariance.

\begin{theore}[\!]
Let $A$ be as in the statement of Theorem~{\rm 1}. Assume that $A$
is hermitian and a $Q$-algebra. Then
$RK_{*}(S(\R,A^{\infty},\alpha)\cong
K_{*}(C^{*}(\R,E(A),\alpha))$. Further if the action $\alpha$ is
isometric on $A${\rm ,} then $RK_{*}(L^{1}(\R,A,\alpha))\cong
K_{*}(C^{*}(\R,E(A),\alpha))$.
\end{theore}

We apply this to the differential structure of a
$C^{*\!}$-algebra. Let $\alpha$ be an action of $\R$ on a
$C^{*\!}$-algebra $A$ leaving a dense {$^{*\!}$}-subalgebra $\cU$
invariant. Let $T\sim (T_{k})_{0}^{\infty}$ be a differential
\hbox{$^{*\!}$-seminorm} on $\cU$ in the sense of Blackadar and
Cuntz \cite{5} with $T_{0}(x) = \|\cdot\|$ the $C^{*\!}$-norm from
$A$. Let $T$ be $\alpha$-invariant. Let $\cU_{(k)}$ be the
completion of $\cU$ in the submultiplicative \hbox{$^{*\!}$-norm}
$p_{k}(x) = \sum_{i=0}^{k} T_{i}(x)$. The differential Frechet
$^{*\!}$-algebra defined by $T$ is $\cU_{\tau} =
\lim\limits_{\leftarrow} \cU_{(k)}$, the inverse limit of Banach
$^{*\!}$-algebras $\cU_{(k)}$.

Now consider $\tilde{\cU}$ to be the $\alpha$-invariant smooth
envelope of \ $\cU$ defined to be the completion of $\cU$ in the
collection of all $\alpha$-invariant differential
$^{*\!}$-seminorms. Notice that neither $\cU_{\tau}$ nor
$\tilde{\cU}$ is a subalgebra of $A$, though each admits a
continuous surjective $^{*\!}$-homomorphism onto $A$ induced by
the inclusion $\cU\rightarrow A$. There exists actions of $\R$ on
each of $\cU_{\tau}$ and $\tilde{\cU}$ induced by $\alpha$. The
following is a smooth Frechet analogue of Connes' analogue of Thom
isomorphism \cite{7}. It supplements an analogues result in
\cite{11}.

\begin{theore}[\!]$\left.\right.$
\begin{enumerate}
\renewcommand{\labelenumi}{\rm (\alph{enumi})}
\leftskip .15pc
\item $RK_{*}(S(\R,\cU_{\tau}^{\infty},\alpha)) = K_{*+1}(A)$.

\item Assume that $\tilde{\cU}$ is metrizable. Then
$RK_{*}(S(\R,\tilde{\cU},\alpha)) = K_{*+1}(A)$.\vspace{-1.5pc}
\end{enumerate}
\end{theore}

\section{Preliminaries and notations}

A {\it Frechet $^{*\!}$-algebra} $(A,t)$ is a complete topological
involutive algebra $A$ whose topology $t$ is defined by a
separating sequence $\{\|\cdot\|_{n}\hbox{:}\ n\in \N\}$ of
seminorms satisfying $\|xy\|_{n}\leq \|x\|_{n}\|y\|_{n},\
\|x^{*}\|_{n} = \|x\|_{n},\ \|x\|_{n}\leq \|x\|_{n+1}$ for all
$x,y$ in $A$ and all $n$ in $\N$. If each $\|\cdot\|_{n}$
satisfies $\|x^{*}x\|_{n} = \|x\|_{n}^{2}$ for all $x$ in $A$,
then $A$ is a $\sigma$-$C^{*\!}$-{\it algebra} \cite{9}. $A$ is
called a \hbox{$Q$-{\it algebra}} if the set of all quasi-regular
elements of $A$ is an open set. For each $n$ in $\N$, let $A_{n}$
be the Hausdorff completion of $(A,\|\cdot\|_{n})$. There exists
norm decreasing surjective \hbox{$^{*\!}$-homomorphisms}
$\pi_{n}\hbox{:}\ A_{n+1}\rightarrow A_{n}$,
$\pi_{n}(x+\ker\|\cdot\|_{n+1}) = x + \ker\|\cdot\|_{n}$ for all
$x\in A$. Then the sequence
\begin{align*}
A_{1} \xleftarrow{\quad\pi_{1}\quad} A_{2} \xleftarrow{\ \
\pi_{2}\ \ } A_{3} \xleftarrow{\ \ \pi_{3}\ \ } \cdots
\xleftarrow{\ \ \pi_{n-1}\ \ } A_{n} \xleftarrow{\ \ \pi_{n}\ \ }
A_{n+1} \longleftarrow \cdots
\end{align*}
is an inverse limit sequence of Banach $^{*\!}$-algebras and
$A=\lim\limits_{\leftarrow} A_{n}$, the inverse limit of Banach
$^{*\!}$-algebras. Let $\Rep(A)$ be the set of all
$^{*\!}$-homomorphisms $\pi\hbox{:}\ A\rightarrow B(H_{\pi})$ of
$A$ into the $C^{*\!}$-algebras $B(H_{\pi})$ of all bounded linear
operators on Hilbert spaces $H_{\pi}$. Let
\begin{align*}
\Rep_{n}(A) &:= \{\pi\in \Rep(A): \ \hbox{there exists}\ k>0\ \hbox{such
that}\\[.5pc]
&\quad\ \ \|\pi(x)\| \leq k \| x\|_{n}\ \hbox{for all}\ x\}.
\end{align*}
Then $|x|_{n} := \sup \{\|\pi(x)\|\hbox{:}\ \pi \in \Rep_{n}(A)\}$ defines a
$C^{*\!}$-seminorm on $A$. The {\it star radical} of $A$ is
\begin{equation*}
\srad(A) = \{x\in A\!:|x|_{n} = 0\ \hbox{for all}\ n \ \hbox{in}\
\N\}.
\end{equation*}
The enveloping $\sigma$-$C^{*\!}$-algebra $(E(A),\tau)$ of $A$ is the
completion of $A/\srad(A)$ in the topology $\tau$ defined by the
$C^{*\!}$-seminorms $\{|\cdot|_{n}\hbox{:}\ n\in \N\},\ |x+\srad(A)|_{n} = |x|_{n}$
for $x$ in $A$.

Let $\alpha$ be a strongly continuous action of $\R$ by continuous
$^{*}$-automorphisms of $A$. The $C^{\infty}$-{\it elements of A for the
action} $\alpha$ are
\begin{equation*}
A^{\infty} := \{x\in A\!: t \rightarrow \alpha_{t}(x)\ \hbox{is
a}\ C^{\infty}\hbox{-function}\}.
\end{equation*}
It is a dense $^{*\!}$-subalgebra of $A$ which is a Frechet algebra with
the topology defined by the submultiplicative $^{*\!}$-seminorms
\begin{equation*}
\|x\|_{k,n} = \|x\|_{n} + \sum\limits_{j=0}^{k} (1/j!) \|
\delta^{j} x\|_{n},\quad n\in\N,\ k\in\Z^{+} = \N\cup (0)
\end{equation*}
where $\delta$ is the derivation $\delta(x) = (\d/\d
t)\alpha_{t}(x)|_{t=0}$. By Theorem~A.2 of \cite{14}, $\alpha$ leaves
$A^{\infty}$ invariant and each $\alpha_{t}$ restricted to $A^{\infty}$
gives a continuous $^{*\!}$-automorphism of the Frechet algebra
$A^{\infty}$. The action $\alpha$ is {\it smooth} if $A^{\infty} = A$.\vspace{.2pc}

\subsection{\it Smooth Schwartz crossed product {\rm \cite{14}}}\vspace{-.4pc}

Assume that $\alpha$ is $m$-{\it tempered} in the sense that for
each $n\in\N$, there exists a polynomial $P_{n}$ such that
$\|\alpha_{r}(x)\|_{n} \leq P_{n}(r) \|x\|_{n}$ for all $r\in \R$
and all $x\in A$. Let $S(\R)$ be the Schwartz space. The completed
(projective) tensor product $S(\R)\otimes A = S(\R,A)$ consisting
of \hbox{$A$-valued} Schwartz functions on $\R$ is a Frechet
algebra with the twisted convolution
\begin{equation*}
(f*g)(r) = \int_{R} f(s)\alpha_{s} (g(r-s))\d s
\end{equation*}
called the {\it smooth Schwartz crossed product by} $\R$ denoted
by $S(\R,A,\alpha)$. The algebra $S(\R,A^{\infty},\alpha)$ is a
Frechet $^{*\!}$-algebra with the involution $f^{*}(r) =
\alpha_{r}(f(-r)^{*})$ (Corollary~4.9 of \cite{14}) whose topology
$\tau_{s}$ is defined by the seminorms
\begin{align*}
\|f\|_{n,l,m} = \sum\limits_{i+j=n} \int_{R} (1+|r|)^{i}\|
f^{(j)}(r) \|_{l,m}\d r,\quad n\in\Z^{+},l\in\Z^{+}, m\in \N
\end{align*}
where
\begin{equation*}
\|f^{(j)} (r)\|_{l,m} = \sum\limits_{k=0}^{l} (1/k!)\|\delta^{k}
(\alpha_{s}((\d^{j}/\d r^{j}) f(r))|_{s=0}\|_{m}
\end{equation*}
(Theorem 3.1.7 of \cite{14}, \cite{11}). These seminorms are
submultiplicative if $\alpha$ is isometric on $A$ in the sense
that $\|\alpha_{r}(x)\|_{n} = \|x\|_{n}$ for all $n\in\N$ and all
$x\in A$.\vspace{.2pc}

\subsection{\it $L^{1}$-crossed products {\rm \cite{11,14}}}\vspace{-.4pc}

Let $F_{d}$ be the set of all functions $f\hbox{:}\ \R\rightarrow A$ for which
\begin{equation*}
\|f\|_{d,m} := \int_{R}(1+|r|)^{d} \|f(r)\|_{m}\d r <\infty
\end{equation*}
for all $m$ in $\N$. Here $\int$ denotes the upper integral. Let
$\lL_{d}$ be the closure in $F_{d}$ of the set of all measurable
simple functions $f\hbox{:}\ \R\rightarrow A$ in the topology on
$F_{d}$ given by the seminorms $\{\|\cdot\|_{d,m}\hbox{:}\ m\in
\N\}$. Let $N_{d}=\cap\{\ker\|\cdot\|_{d,m}\hbox{:}\ m\in\N\}$.
Then $N_{d}=N_{d+1}$; $L_{d}:=\lL_{d}/N_{d}$ is complete in
$\{\|\cdot\|_{d,m}\hbox{:}\ m\in\N\}$ and $L_{d+1}\rightarrow
L_{d}$ continuously. The {\it space of $|\cdot|$-rapidly vanishing
$L^{1}$-functions from $\R$ to} $A$ is
$L^{1}_{|\cdot|}(\R,A,\alpha):= \cap\{L_{d}\hbox{:}\
d\in\Z^{+}\}$, a Frechet algebra with the topology given by the
seminorms $\{\|\cdot\|_{d,m}\hbox{:}\ m\in\N,\ d\in\Z^{+}\}$ and
with twisted convolution. Assume that $\alpha$ is isometric on
$(A,\{\|\cdot\|_{n}\})$. Then the completed projective tensor
product $L^{1}(\R)\otimes A = L^{1}(\R,A)$ is a Frechet
$^{*\!}$-algebra with twisted convolution and the involution
$f\rightarrow f^{*}$. This $L^{1}$-{\it crossed product} is
denoted by $L^{1}(\R,A,\alpha)$. Notice that $\alpha$ is isometric
on $(A^{\infty},\{\|\cdot\|_{n,m}\})$ also, so that the Frechet
$^{*\!}$- algebra $L^{1}(\R,A^{\infty},\alpha)$ is defined; and
then the induced actions $(\alpha_{r}f)(s) = \alpha_{r}(f(s))$ on
$L^{1}(\R,A^{\infty},\alpha)$ and on $L^{1}(\R,A,\alpha)$ are also
isometric.\vspace{.4pc}

\subsection{\it $\sigma$-$C^{*\!}$-crossed product}\vspace{-.3pc}

Assume that $\alpha$ is isometric. We define the
$\sigma$-$C^{*\!}$-{\it crossed product} $C^{*}(\R,A,\alpha)$ of
$A$ by $\R$ to be the enveloping $\sigma$-$C^{*\!}$-algebra
$E(L^{1}(\R,A,\alpha))$ of $L^{1}(\R,A,\alpha)$.

\section{Technical lemmas}

\begin{lem}
Let $\alpha$ be $m$-tempered on $A$. Then $\alpha$ extends as a strongly
continuous isometric action of $\R$ by continuous $^{*\!}$-automorphisms
of the $\sigma$-$C^{*\!}$-algebra $E(A)$.
\end{lem}

\begin{proof}
By the $m$-temperedness of $\alpha$, for each $n\in\N$, there exists a
polynomial $P_{n}$ such that for all $x\in A$ and all $r\in\R$,
$\|\alpha_{r}(x)\|_{n}\leq P_{n}(r)\|x\|_{n}$. Let $r\in\R$. Let
$x\in\srad(A)$. Then for all $\pi\in\Rep(A),\ \pi(x) = 0$, so that
$\sigma(\alpha_{r}(x)) = 0$ for all $\sigma\in\Rep(A)$, hence
$\alpha_{r}(x)\in \srad(A)$. Thus $\alpha_{r}(\srad(A))\subseteq
\srad(A)$, and the map
\begin{equation*}
\tilde{\alpha}_{r}\!: A/\srad(A) \rightarrow A /\srad(A),\quad
\tilde{\alpha}_{r}([x]) = [\alpha_{r}(x)],
\end{equation*}
where $[x]=x+\srad(A)$, is a well-defined $^{*\!}$-homomorphism.
Further, let $\tilde{\alpha}_{r}[x]=0$. Then
$\alpha_{r}(x)\in\srad(A)$. Hence
$x=\alpha_{-r}(\alpha_{r}(x))\in\srad(A),\ [x] = 0$. Thus
$\tilde{\alpha}_{r}$ is one-to-one, which is clearly surjective.
Now, for each $n\in\N$, and for all $x\in A$,
\begin{equation*}
|\tilde{\alpha}_{r}[x]|_{n} = |[\alpha_{r}(x)]|_{n}\leq
\|\alpha_{r}(x)\|_{n} \leq P_{n}(r)\|x\|_{n}.
\end{equation*}
Since, by definition, $|\cdot|_{n}$ is the greatest $C^{*\!}$-seminorm on
$A/\srad(A)$ satisfying that for some $k_{n}>0,\ |[z]|_{n}\leq
k_{n}\|z\|_{n}$ for all $z\in A$, it follows that
$|\tilde{\alpha}_{r}[x]|_{n}\leq |[x]|_{n}$ for all $x$ in $A$. Hence
\begin{equation*}
|[x]|_{n} \leq |\tilde{\alpha}_{-r} (\tilde{\alpha}_{r}[x]|_{n} =
|\tilde{\alpha}_{-r} [\alpha_{r}(x)]|_{n}\leq |[\alpha_{r}(x)]|_{n} = |\tilde{\alpha}_{r}[x]|_{n}
\end{equation*}
showing that $|\tilde{\alpha}_{r}[x]|_{n} = |[x]|_{n}$ for all
$x\in A,\ r\in\R,n\in\N$. It follows that $\tilde{\alpha}_{r}$
extends as a $^{*\!}$-automorphism $\tilde{\alpha}_{r}\hbox{:}\
E(A)\rightarrow E(A)$ satisfying $|\tilde{\alpha}_{r}(z)|_{n} =
|z|_{n}$ for all $z\in A$ and all $n\in\N$; and
$\tilde{\alpha}\hbox{:}\ \R\rightarrow \hbox{Aut}^{*}(E(A)),\
r\rightarrow \tilde{\alpha}_{r}$ defines an isometric action of
$\R$ on $E(A)$. We verify that $\tilde{\alpha}$ is strongly
continuous. Let $z\in E(A)$. It is sufficient to prove that the
map $f\hbox{:}\ \R\rightarrow E(A),\ f(r) = \alpha_{r}(z)$ is
continuous at $r=0$. Choose $z_{n} = [x_{n}]$ in $A/\srad(A)$ such
that $z_{n}\rightarrow z$ in $E(A)$. Fix $k\in\N,\ \varepsilon
>0$. Choose $n_{0}$ in $\N$ such that $|z_{n_{0}} - z|_{k} <
\varepsilon/3$ with $z_{n_{0}} = [x_{n_{0}}]$. Then for all
$r\in\R$, $|\tilde{\alpha}_{r}(z) -
\tilde{\alpha}_{r}(z_{n_{0}})|_{k} = |z-z_{n_{0}}|_{k} <
\varepsilon/3$. Since $\alpha$ is strongly continuous, there
exists a $\delta>0$ such that $|r|<\delta$ implies that
$\|\alpha_{r}(x_{0}) - x_{0}\|_{k} <\varepsilon /3$. Then for all
such $r$, $|\tilde{\alpha}_{r}(z)-z|_{k}<\varepsilon$ showing the
desired continuity of $f$. This completes the proof.\hfill
$\Box$\vspace{.7pc}
\end{proof}

\begin{notat}
Henceforth we denote the action $\tilde{\alpha}$ by $\alpha$.

A {\it covariant representation} of the Frechet algebra dynamical system
$(\R,A,\alpha)$ is a triple $(\pi,U,H)$ such that

\begin{enumerate}
\renewcommand{\labelenumi}{(\alph{enumi})}
\leftskip .15pc
\item $\pi\hbox{:}\ A\rightarrow B(H)$ is a $^{*\!}$-homomorphism;
\item $U\hbox{:}\ \R\rightarrow \cU(H)$ is a strongly continuous unitary
representation of $\R$ on $H$; and
\item $\pi(\alpha_{t}(x)) = U_{t}\pi(x)U_{t}^{*}$ for all $x\in A$ and all
$t\in\R$.\vspace{-.5pc}
\end{enumerate}
\end{notat}

The following is an analogue of Proposition~7.6.4, p.~257 of
\cite{12} which can be proved along the same lines. Let
$C_{c}^{\infty}(\R,A^{\infty}) = C_{c}^{\infty}(\R)\otimes
A^{\infty}$ (completed projective tensor product) be the space of
all $A^{\infty}$-valued $C^{\infty}$-functions on $\R$ with
compact supports.

\begin{lem}
Let $A$ have a bounded approximate identity $(e_{l})$ contained in
$A^{\infty}$ which is also a bounded approximate identity for the
Frechet algebra $A^{\infty}$. {\rm (}In particular{\rm ,} let $A$ be unital.{\rm )}

\begin{enumerate}
\renewcommand{\labelenumi}{\rm (\alph{enumi})}
\leftskip .15pc
\item If $(\pi,U,H)$ is a covariant representation of
$(\R,A^{\infty},\alpha)${\rm ,} then there exists a non-degenerate $^{*\!}$-
representation $(\pi\times U,H)$ of $S(\R,A^{\infty},\alpha)$ such that
\begin{equation*}
\hskip -1.25pc (\pi\times U)y = \int_{R}\pi(y(t)) U_{t}\d t
\end{equation*}
for every $y$ in $C_{c}^{\infty}(\R,A^{\infty})$. The correspondence
$(\pi,U,H)\rightarrow (\pi\times U,H)$ is bijective onto the set of all
non-degenerate $^{*\!}$-representations of $S(\R,A^{\infty},\alpha)$.

\item Let $\alpha$ be isometric. Then the above gives a one-to-one
correspondence between the covariant representations of $(\R,A,\alpha)$
and non-degenerate $^{*\!}$-representations of each of
$L^{1}(\R,A^{\infty},\alpha)$ and $L^{1}(\R,A,\alpha)$.
\end{enumerate}
\end{lem}

\begin{lem}
$E(A^{\infty}) = E(A)${\rm ;} and for all $k$ in $\Z^{+},n$ in $\N,
\|_{n,k} = \|_{n}$.
\end{lem}

\begin{proof}
Consider the inverse limit $A = \lim\limits_{\leftarrow} A_{n}$ as
in the Introduction. Since $\alpha$ satisfies
$\|\alpha_{r}(x)\|_{n} \leq P_{n}(r)\|x\|_{n}$ for all $x\in\R$,
all $r\in A$ and all $n\in\N$, it follows that for each $n,\alpha$
`extends' uniquely as a strongly continuous action $\alpha^{(n)}$
of $\R$ by continuous \hbox{$^{*\!}$-automorphisms} of the Banach
$^{*\!}$-algebra $A_{n}$. Let $(A_{n,m}, \|\cdot\|_{n,m})$ be the
Banach algebra consisting of all $C^{m}$-elements $y$ of $A_{n}$
with the norm $\|\cdot\|_{n,m} = \|y\|_{n} + \sum_{i=1}^{m} (1/i!)
\|\delta^{i}(x)\|_{n}$. Let $(A_{n}^{\infty},\ \{\|\
\|_{m,n}\hbox{:}\ m\in\Z^{+}\})$ be the Frechet algebra consisting
of all $C^{\infty}$-elements of $A_{n}$ for the action
$\alpha^{(n)}$. Then
\begin{equation*}
A^{\infty} = \lim\limits_{\leftarrow} A_{n}^{\infty} =
\lim\limits_{\leftarrow}\lim\limits_{\leftarrow} A_{m,n} =
\lim\limits_{\leftarrow} A_{n,n}.
\end{equation*}
By Theorem~2.2 of \cite{15}, each $A_{m,n}$ is dense and spectrally
invariant in $A_{n}$. Hence each $A_{n,m}$ is a $Q$-normed algebra in
the norm $\|\cdot\|_{n}$ of $A_{n}$.

Let $\pi\hbox{:}\ A^{\infty}\rightarrow B(H)$ be a $^{*\!}$-representation of $A$
on a Hilbert space $H$. Since the topology of $A^{\infty}$ is determined
by the seminorms
\begin{equation*}
\|x\|_{n,n} = \|x\|_{n} + \sum\limits_{j=1}^{n} (1/j!)\|\delta^{j}
(x)\|_{n},\quad n\in\N
\end{equation*}
it follows that for some $k>0,\ \|\pi(x)\|\leq k\|x\|_{n,n}$ for all
$a\in A^{\infty}$. Hence $\pi$ defines a $^{*\!}$-homomorphism
$\pi\hbox{:}\ (A_{n,n}, \|\ \|_{n,n})\rightarrow B(H)$ satisfying $\|\pi(x)\|\leq
k\|x\|_{n,n}$ for all $x$ in $A_{n,n}$. Since $(A_{n,n}, \|\ \|_{n})$ is
a $Q$-normed $^{*\!}$-algebra, this map $\pi$ is continuous in the norm
$\|\ \|_{n}$ on $A_{n,n}$. In fact, for all $x$ in $A^{\infty}$,
\begin{align*}
\|\pi(x)\|^{2} &= \|\pi(x^{*} x)\| = r_{B(H)} (\pi (x^{*} x)) \leq
r_{A_{n,n}} (\pi(x^{*}x + \ker\|\ \|_{n,n}))\\[.3pc]
&\leq \|x^{*} x+ \ker \|\ \|_{n} \|=\|x^{*} x\|_{n}\leq \|x\|^{2}.
\end{align*}
Thus $\|\pi(x)\|\leq \|x\|_{n}$ for all $x$ in $A^{\infty}$. Since
$A^{\infty}$ is dense in $A$, $\pi$ can be uniquely extended as a
$^{*\!}$-representation $\pi\hbox{:}\ A\rightarrow B(H)$ satisfying that
$\|\pi(x)\|\leq \|x\|_{n}$ for all $x$ in $A$. Then by the definition of
the $C^{*\!}$-seminorm $|\ |_{n}$ on $A,\pi$ extends as a continuous
$^{*\!}$-homomorphism $\tilde{\pi}\hbox{:}\ E(A)\rightarrow B(H)$ such that
$\|\tilde{\pi}(x)\|\leq |x|_{n}$ for all $x$ in $E(A)$. This also
implies that $E(A^{\infty}) = E(A)$ and $|\cdot|_{n,m} = |\cdot|_{n}$
for all $n,m$.
\end{proof}

\begin{lem}
Let $B$ be a $\sigma$-$C^{*\!}$-algebra. Let $j\hbox{\rm :}\ A\rightarrow E(A)$ be
$j(x) = x+\srad(A)$. Let $\pi\hbox{\rm :}\ A\rightarrow B$ be a $^{*\!}$-homomorphism.
Then there exists a unique $^{*\!}$-homomorphism $\tilde{\pi}\hbox{\rm :}\
E(A)\rightarrow B$ such that $\pi = \tilde{\pi} \circ j$.
\end{lem}

This follows immediately by taking
$B=\lim\limits_{\leftarrow}B_{n}$, where $B_{n}$'s are
$C^{*\!}$-algebras, and by the universal property of
$E(A)$.\vspace{-.5pc}

\section{Proof of Theorem 1}

\begin{step}
\renewcommand\thestep{{\it \arabic{Step}}}
{\rm $\Rep(S(\R, A^{\infty},\alpha)) = \Rep(S(\R,E(A),\alpha))=\Rep
(L^{1}(\R, E(A),\alpha))$ up to one-to-one correspondence.}
\end{step}

By Lemma~3.1, the Frechet algebras $S(\R,E(A),\alpha)$ and
$L^{1}(\R,E(A),\alpha)$ are \hbox{$^{*\!}$-algebras} with the
continuous involution $y\rightarrow y^{*}, y^{*}(t) =
\alpha_{t}(y(-t))^{*}$. By Lemma~3.2, $\Rep(S(\R, E(A), \alpha)) =
\Rep(L^{1}(\R, E(A),\alpha))$ each identified with the set of all
covariant representations. Let $\rho\hbox{:}\ S(\R,
A^{\infty},\alpha)\rightarrow B(H)$ be in
$\Rep(S(\R,A^{\infty},\alpha))$. There exists $c>0$ and
appropriate $n,l,m$ such that for all $y$,
\begin{equation}
\|\rho(y)\| \leq c\|y\|_{n,l,m} = c\sum\limits_{i+j=n}\int_{R}
(1+|r|)^{i}\| y^{(j)}(r)\|_{l,m}\d r.
\end{equation}
By Lemma~3.2, there exists a covariant representation $(\pi,U,H)$
of $(\R,A^{\infty},\alpha)$ on $H$ such that $\rho=\pi\times U$.
Thus $\pi\hbox{:}\ A^{\infty}\rightarrow B(H)$ is a
$^{*\!}$-homomorphism and $U\hbox{:}\ \R\rightarrow \cU(H)$ is a
strongly continuous unitary representation such that

\begin{enumerate}
\renewcommand{\labelenumi}{(\roman{enumi})}
\leftskip .4pc
\item  \hskip 1cm $\rho(f) = \int_{R}\pi (f(t)) U_{t}\d t\quad
\hbox{for all} \ f\ \hbox{in}\ S(\R,A^{\infty},\alpha),$\hfill
(2)\vspace{.5pc}

\item \hskip 1cm $\pi(\alpha_{t}(x)) = U_{t}\pi(x)U_{t}^{*}\quad\hbox{for all} \ x\in
A^{\infty}, t\in \R,$\hfill (3)\vspace{.5pc}

\item  \hskip 1cm there exists $K>0$ such that $\|\pi(x)\| \leq k\|x\|_{l,m}$ for
all $x\in A^{\infty}$.
\end{enumerate}

\noindent The $l,m$ in (iii) are the same as in (1). Let
$\{|\cdot|_{l,m}\hbox{:}\ l$ in $\Z^{+},m$ in $\N\}$ be the
sequence of \hbox{$C^{*\!}$-seminorms} on $A^{\infty}$ (and also
on $E(A^{\infty})$ via $\srad\,A^{\infty}$) which are defined by
the submultiplicative $^{*\!}$-seminorms
$\{\|\cdot\|_{l,m}\hbox{:}\ l$ in $\Z^{+},m$ in $\N\}$. Then
$|\cdot|_{l,m}$ is the greatest \hbox{$C^{*\!}$-seminorm} on
$A^{\infty}$ satisfying that there exists $M=M_{l,m} >0$ such that
$|\cdot|_{l,m}\leq$ \hbox{$M\|\cdot\|_{l,m}$}. Hence by (iii)
above, $\pi$ can be uniquely extended as a continuous
\hbox{$^{*\!}$-homomorphism} $\tilde{\pi}\hbox{:}\
E(A^{\infty})\rightarrow B(H)$ such that $\tilde{\pi}(j(x)) =
\pi(x)$ for all $x\in A^{\infty}$; and \setcounter{equation}{3}
\begin{equation}
\hskip -1.25pc  \|\tilde{\pi}(x)\| \leq |x|_{l,m}\ \ \hbox{for all}\ x\in E(A^{\infty}).
\end{equation}
Here $j$ is the map $j\hbox{:}\ A^{\infty}\rightarrow
E(A^{\infty})$, $j(x)=x+\srad\, A^{\infty}$. Let $l$ denote
$\max(l,m)$. Then we have
\begin{align}
\hskip -1.25pc  \|\rho(y)\| &\leq c \|y\|_{n,l,l}\ \ \hbox{for all}\ y\in
S(\R,A^{\infty},\alpha);\nonumber\\[.2pc]
\hskip -1.25pc  \|\pi(x)\| &\leq k \|x\|_{l,l}\ \ \hbox{for all}\ x\in A^{\infty};\nonumber\\[.2pc]
\hskip -1.25pc  \|\tilde{\pi}(z)\| &\leq |z|_{l,l}\ \ \hbox{for all}\ z \in
E(A^{\infty}).
\end{align}
By Lemma~3.3, $\tilde{\pi}\hbox{:}\ E(A)\rightarrow B(H)$ is a
$^{*\!}$-representation satisfying $\|\tilde{\pi}(x)\|\leq
|x|_{l}$ for all $x$ in $E(A)$. We have the following commutative
diagram.\vspace{.7pc}

\hskip 2.5pc{\epsfbox{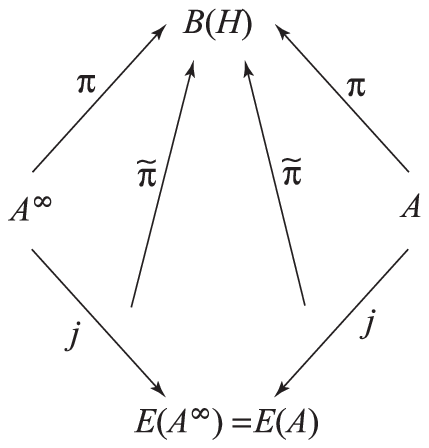}}\vspace{.7pc}

Now, let $\alpha\hbox{:}\ \R\rightarrow\hbox{Aut}^{*} E(A)$ be the action on
$E(A)$ induced by $\alpha$ as in Lemma~3.1 satisfying
\begin{equation}
\hskip -1.25pc  \alpha_{t}(j(x)) = j(\alpha_{t}(x))\quad \hbox{for all}\ x\ \hbox{in}\
A.
\end{equation}
Then $(\tilde{\pi},U,H)$ is a covariant representation of
$(\R,E(A),\alpha)$. Indeed, let $x\in A^{\infty}, y = j(x)$. Then for
all $t\in\R$,
\begin{align*}
\hskip -1.25pc \tilde{\pi}(\alpha_{t}(y)) &= \tilde{\pi}(\alpha_{t}(j(x))) =
\tilde{\alpha}(j(\alpha_{t}(x))) = \pi(\alpha_{t}(x)) =
U_{t}\pi(x)U_{t}^{*}\\[.2pc]
&= U_{t}\tilde{\pi} (j(x))U_{t}^{*} = U_{t}\tilde{\pi}(y)U_{t}^{*}.
\end{align*}
By the continuity of $\tilde{\pi}$ and $\alpha_{t}$, it follows
that $\tilde{\pi}(\alpha_{t}(y)) = U_{t}\tilde{\pi}(y)U_{t}^{*}$
for all $y\in E(A)$ and all $t\in\R$. Hence by Lemma~3.2,
$\tilde{\rho} = \tilde{\pi}\times U$ is a non-degenerate
$^{*\!}$-representation of each of $S(\R,E(A),\alpha)$ and
$L^{1}(\R,E(A),\alpha)$ satisfying, for some constants $c$ and
$c'$, the following (using (5)):

\begin{enumerate}
\renewcommand\labelenumi{(\roman{enumi})}
\setcounter{enumi}{3}
\leftskip .4pc
\item \hskip 1cm For all $f$ in $L^{1}(\R, E(A),\alpha)$, $\|\rho(f)\|\leq c|f|_{l}
= c\int_{R}|f(t)|_{l} \d t$.\vspace{.3pc}

\item \hskip 1cm For all $f\ \hbox{in}\ S(\R, E(A), \alpha),
\|\tilde{\rho}(f)\|\leq c'|f|_{n,l,m}$.\hfill (7)
\end{enumerate}
Thus given a $^{*\!}$-representation $\rho$ of $S(\R,A^{\infty},\alpha)$,
there is canonically associated a \hbox{$^{*\!}$-representation} $\tilde{\rho}$
of each of $S(\R,E(A),\alpha)$ and $L^{1}(\R,E(A),\alpha)$.

Conversely, given $\rho$ in $\Rep(S(\R,E(A),\alpha))$, $\rho = \pi\times
U$ for a covariant representation $(\pi,U)$ of $(\R,E(A),\alpha)$, $\pi
\circ j$ is a covariant representation of $A$, and then $(\pi\circ
j)\times U$ is in $\Rep(S(\R,A^{\infty},\alpha))$.

\begin{step}
{\rm The $\sigma$-$C^{*\!}$-algebra $C^{*}(\R,E(A),\alpha)$ is universal for
the $^{*\!}$-representations of the Frechet algebra
$S(\R,A^{\infty},\alpha)$.}
\end{step}

Let $\tilde{j}\hbox{:}\ S(\R,A^{\infty},\alpha)\rightarrow
L^{1}(\R, E(A),\alpha)$ be the map
\setcounter{equation}{7}
\begin{align}
\tilde{j}(f) &= j\circ f = \tilde{f}\ \ \hbox{(say)}, \quad
\hbox{i.e.,}\nonumber\\[.2pc]
\tilde{j}(f)(r) &= j(f(r)) = f(r) + \srad(A)^{\infty} \quad
\hbox{for all}\ \ r\in \R.
\end{align}
Notice that the map $\tilde{j}$ is defined and is continuous;
because $(S(\R,A^{\infty},\alpha))\subset L^{1}(\R,$
$A^{\infty},\alpha)\subset L^{1}(\R,A,\alpha)$, and for $n$ in
$\N$ and $m$ in $\Z^{+}$, all $f$ in $S(\R,A^{\infty},\alpha)$,
\begin{align*}
&|\tilde{f}(t)|_{n} \leq \|f(t)\|_{n} \leq M\|f(t)\|_{m,n},\ \
\hbox{and
hence}\\[.2pc]
&\int_{R}|\tilde{f}(t)|_{l}\d t \leq \int_{R} \|f(t)\|_{m,n}\d
t < \infty
\end{align*}
so that $f\in L^1(\R, E(A), \alpha)$. Let $j_1\hbox{:}\ L^1(\R,
E(A),\alpha)\rightarrow C^* (\R, E(A),\alpha)$ be the natural map
$j_1(f) = f + \hbox{srad}(L^1(\R, E(A),\alpha))$. This gives the
continuous $^*$-homomorphism
\begin{equation}
J : j_1 \circ \tilde{j} : S(\R, A^\infty, \alpha) \rightarrow C^*
(\R, E(A), \alpha).
\end{equation}

\hskip 3.3pc {\epsfbox{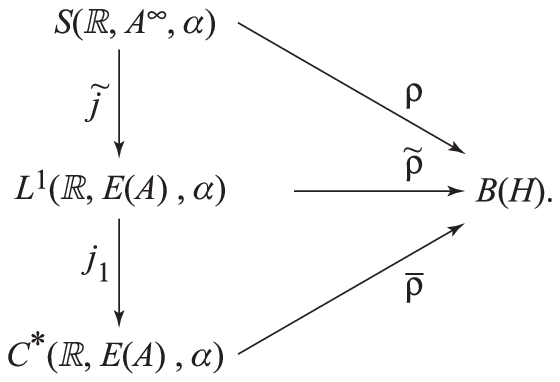}}\vspace{.7pc}

Let $\rho \in \hbox{Rep}(S(\R, A^\infty, \alpha)), \rho = \pi
\times U$ in usual notations with $\pi\hbox{:}\ A^\infty
\rightarrow B(H)$ in $\hbox{Rep}(E(A))$ such that $\pi=\tilde{\pi}
\circ j$. Let $\tilde{\rho}\hbox{:}\ L^1(\R, E(A), \alpha)
\rightarrow B(H)$ be $\tilde{\rho} = \tilde{\pi} \times U$. Then
for all $f$ in $S(\R, A^\infty, \alpha)$,
\begin{align*}
\tilde{\rho} (\tilde{j}(f)) &= (\tilde{\pi} \times U)
(\tilde{j}(f)) = \int_R \tilde{\pi} (\tilde{j}(f)(t)) U_t \d t =
\int_R \tilde{\pi}(j \circ f) (t) U_t \d t\\[.5pc]
&= \int_R \tilde{\pi} (j(f(t))) U_t \d t = \int_R \tilde{\pi}
(f(t) + \hbox{srad}(A)) U_t \d t\\[.5pc]
&= \int_R \pi(f(t)) U_t \d t = \rho (f).
\end{align*}
Thus $\tilde{j}\circ \tilde{\rho} = \rho$; and hence $J \circ
\bar{\rho} = \rho$, where $J = j_1 \circ \tilde{j}$ and
$\bar{\rho} \in \hbox{Rep}(C^*(\R, E(A), \alpha))$ is defined by
$j_1 \circ \bar{\rho} = \tilde{\rho}$ in view of $C^*(\R, E(A),
\alpha) = E(L^1(\R, E(A), \alpha))$.

\begin{step}
{\rm Given a $^*$-homomorphism $\rho\hbox{:}\ S(\R, A^{\infty},
\alpha) \rightarrow B$ from $S(\R, A^{\infty}, \alpha)$ to a
\hbox{$\sigma$-$C^*$-algebra} $B$, there exists $^*$-homomorphisms
$\tilde{\rho}\hbox{:}\ L^1(\R, E(A),\alpha) \rightarrow B,
\tilde{\rho}\hbox{:}\ C^*(\R, E(A),\break \alpha) \rightarrow B$ such
that $\rho = \tilde{\rho} \circ \tilde{j} = \bar{\rho} \circ J$
and $\tilde{\rho} = \bar{\rho} \circ j_1$.}
\end{step}

This follows by applying Step~II to each of the factor
$C^*$-algebra $B_n$ in the inverse limit decomposition of $B$.

\begin{step}
{\rm $C^*(\R, E(A), \alpha) = E(S(\R, A^\infty, \alpha))$ up to
homeomorphic $^*$-isomorphism.}
\end{step}
\pagebreak

Let $k\hbox{:}\ S(\R, E(A), \alpha)\rightarrow E(S(\R, A^\infty,
\alpha))$ be $k(f) = f + \hbox{srad}\ S(\R, A^{\infty}, \alpha)$.
Then there exists a $^*$-homomorphism $\bar{k}\hbox{:}\ C^*(\R,
E(A), \alpha)\rightarrow E(S(\R, A^\infty, \alpha))$ such that
$\bar{k} \circ J = k$. We show that $\bar{k}$ is the desired
homeomorphic $^*$-isomorphism making the following diagram
commutative.
\begin{equation*}
{\epsfbox{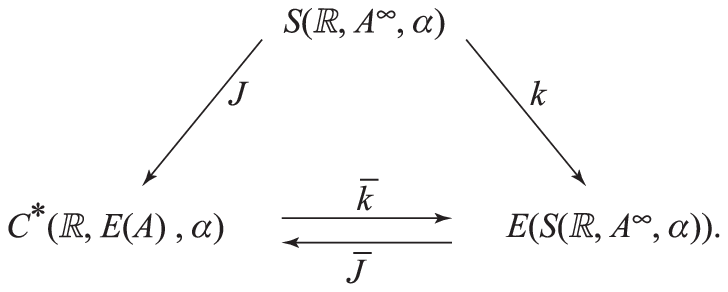}}
\end{equation*}\vspace{-1.5cm}

\hfill (10)\vspace{.6cm}

By the universal property of $E(S(\R, A^\infty, \alpha))$, there
exists a $^*$-homomorphism\break $\bar{J}\hbox{:}\
E(S(\R,A^\infty,\alpha)) \rightarrow C^*(\R, E(A),\alpha)$ such
that $\bar{J} \circ k = J$. We claim that $|\bar{k}|_{{\rm
Im}(J)}$ is injective. Indeed, let $f\in S(\R,A^\infty, \alpha)$
be such that $\bar{k}(J(f)) =0$. Hence $k(f) =0$, so that $f\in
\hbox{srad}(S(\R,A^\infty,\alpha))$. Thus, for all
$\rho\in\hbox{Rep}(S(\R,A^\infty,\alpha))$, $\rho(f)=0$.
Therefore, by Step~I, $\sigma(\bar{f}) =0$ for all
$\sigma\in\hbox{Rep}(L^1(\R,E(A),\alpha))$. (Recall that
$\tilde{f}=j \circ f = \tilde{j}(f)$.) Hence $\tilde{j}(f)$ is in
$\hbox{srad}(L^1(\R,E(A),\alpha))$, and so $j_1(\tilde{j}(f)) =0$.
Therefore $J(f) =0$. It follows that $\bar{k}$ is injective on
$\hbox{Im}(J)$.

Now by (10) and the injectivity of $\bar{k}$ on $\hbox{Im}(J),
\bar{J}\circ k = J$. Hence $J = \bar{J} \circ \bar{k} \circ J$,
and so $\bar{J} \circ\bar{k} = \hbox{id}$ on $\hbox{Im}(J)$.
Similarly $\bar{k}\circ\bar{J}(k(f))=\bar{k}(J(f))=k(f)$, hence
$\bar{k}\circ \bar{J} =\hbox{id}$ on $\hbox{Im}(k)$. Thus $\bar{k}
= (\bar{J})^{-1}$ on $\hbox{Im}(J)$. Thus $\bar{k}$ is a
homeomorphic $^*$-isomorphism from the dense $^*$-subalgebra
$J(S(\R, A^\infty, \alpha))$ of $C^*(\R, E(A), \alpha)$ on the
dense $^*$-subalgebra $k(S(\R, A^\infty, \alpha))$ of
$E(S(\R,A^\infty, \alpha))$. It follows that $C^*(\R,
E(A),\alpha)$ is homeomorphically $^*$-isomorphic to $E(S(\R,
A^\infty, \alpha))$.

\begin{step}
{\rm $E(L_{|\cdot|}^1 (\R, A^\infty, \alpha)) = C^*(\R, E(A),
\alpha)$.}
\end{step}

Let $\R$ act on $L_{|\cdot|}^1 (\R, A, \alpha)$ by $xf(y)=f(x-y)$.
For this action, $(L_{|\cdot|}^1 ((\R, A, \alpha))^\infty = S(\R,
A, \alpha)$ by Theorem~2.1.7 of \cite{14}. Thus $S(\R, A^\infty,
\alpha) = (L_{|\cdot|}^1 ((\R, A, \alpha))^\infty$. Hence by
Lemma~3.4, $E(L_{|\cdot|}^1 (\R, A^\infty, \alpha))^\infty =
E(S(\R, A^\infty,\alpha)) = C^* (\R, E(A), \alpha)$. This
completes the proof of Theorem~1.\hfill $\Box$

\section{Proof of Theorem~2}

Let the Frechet algebra $A$ be hermitian and a $Q$-algebra. Hence
$A$ is spectrally bounded, i.e., the spectral radius $r(x) =
r_A(x) < \infty$ for all $x \in A$. Let $s_A(x) :=
r(x^{*}x)^{1/2}$ be the Ptak's spectral function on $A$. By
Corollary~2.2 of \cite{1}, $E(A)$ is a $C^*$-algebra, the complete
$C^*$-norm of $E(A)$ being defined by the greatest $C^*$-seminorm
$p_\infty(\cdot)$ (automatically continuous) on $A$. Now for any
$x\in A$,
\begin{align*}
p_\infty(x)^2 &= p_\infty(x^{*}x) = \|x^{*}x+ \hbox{srad}(A)\|\\
&= r_{E(A)} (x^{*}x + \hbox{srad}(A)) \leq r_A(x^{*}x) = s_A(x)^2.
\end{align*}
Hence $p_\infty(x) \leq s_A(x)$ for all $x\in A$. By the
hermiticity and $Q$-property, $s_A(\cdot)$ is a
\hbox{$C^*$-seminorm} (Theorem~8.17 of \cite{8}), hence
$p_\infty(\cdot)=s(\cdot) \geq r(\cdot)$. In this case,
$\hbox{rad}(A) = \hbox{srad}(A)$. Let $A_{q} = A/\hbox{rad}(A)$
which is a dense $^*$-subalgebra of the $C^*$-algebra $E(A)$ and
is also a Frechet $Q$-algebra with the quotient topology $t_q$.
Let $[x] = x + \hbox{rad}(A)$ for all $x\in A$. Since the spectrum
\begin{equation*}
\hbox{sp}_A(x) = \hbox{sp}_{A_{q}} ([x]),\quad r_A(x) =
r_{A_q}([x]),\quad s_A(x) = s_{A_q}([x]),
\end{equation*}
and so $r_{A_q}([x])\leq s_{A_q}([x])=\|[x]\|_{\infty}$. Hence
$\|\cdot\|_\infty$ is a spectral norm on $A_q$, i.e.,
$(A_q,\|\cdot\|_\infty)$ is a $Q$-algebra. Thus $A_q$ is
spectrally invariant in $E(A)$. Hence by Corollary~7.9 of
\cite{10}, $K_*(A_q)=RK_*(A_q)=K_*(E(A))$.

Now consider the maps
\begin{equation*}
A \xrightarrow{\ \ \ j \ \ \ } A_{q} \xrightarrow{\ \ \ {\rm id}\
\ \ } E(A)
\end{equation*}
and, for each positive integer $n$, the induced maps
\begin{align*}
M_n(A) \xrightarrow{\ j_n=j\otimes {\rm id}_n} M_n(A_q) =
[M_n(A)]_q \xrightarrow{\ \ \ {\rm id}\ \ \ } M_n(E(A)) =
E(M_n(A)).
\end{align*}
By the spectral invariance of $A$ in $A_q$ via the map $j,
j(\hbox{inv}(A)) = \hbox{inv}(A_q)$, where $\hbox{inv}(K)$ denotes
the group of invertible elements of $K$. Let $\hbox{inv}_0(\cdot)$
denote the principle component in $\hbox{inv}(\cdot)$. We use the
following.

\setcounter{theora}{0}
\begin{lem}
Let $B$ be a Frechet $Q$-algebra or a normed $Q$-algebra. Then
$\hbox{\rm inv}_0(B)$ is the subgroup generated by the range $\exp
B$ of the exponential function.
\end{lem}

The Frechet $Q$-algebra case follows by adapting the proof of the
corresponding Banach algebra result in Theorem~1.4.10 of
\cite{13}. If $(B, \|\cdot\|)$ is a $Q$-normed algebra, then
$(B,\|\cdot\|)$ is advertably complete in the sense that if a
Cauchy sequence $(x_n)$ converges to an element
$x\in\hbox{inv}(B^\sim)$ ($B^\sim$ being the completion of $B$),
then $x\in B$. Hence the exponential function is defined on $B$;
and then the Banach algebra proof can be adapted.

We use the above lemma to verify the following:

\begin{claim}
$j_n(\hbox{inv}_0(M_n(A))) = \hbox{inv}_0(M_n(A_q))$.\vspace{.7pc}
\end{claim}

Take $n=1$. It is clear that $j(\hbox{inv}_0(A)) \subseteq
\hbox{inv}_0(A_q)$. Let $y \in\hbox{inv}_0(A_q)$. Hence $y =\Pi
\exp (z_i)$ for finitely many $z_i=[x_i] = x_i + \hbox{rad}(A)$
for some $x_i$ in $A$. Then $y = [\Pi \exp(x_i)]$. Hence $y\in
j(\hbox{inv}_0(A))$. Thus $j(\hbox{inv}_0(A)) =
\hbox{inv}_0(A_q)$. Now take $n>1$. As $A$ is spectrally invariant
in $A_q$, it follows from Theorem~2.1 of \cite{16} that the
Frechet $Q$-algebra $M_n(A)$ is spectrally invariant in $M_n(A_q)$
via $j_n$. Also, $M_n(A_q)=(M_n(A))_q$ is a $Q$-algebra in both
the quotient topology as well as the $C^*$-norm induced from
$M_n(E(A)) = E(M_nA)$. Applying arguments analogous to above, it
follows that \hbox{$j_n(\hbox{inv}_0(M_n(A)))
=\hbox{inv}_0(M_n(A_q))$.}

Now consider the surjective group homomorphisms
\begin{equation*}
\hbox{inv}(M_n(A)) \xrightarrow{\ \ j_n\ \ }\hbox{inv}(M_n(A_q))
\xrightarrow{\ \ J\ \ }
\hbox{inv}(M_n(A_q))/\hbox{inv}_0(M_n(A_q)).
\end{equation*}
It follows that $\hbox{ker}(J\circ j_{n}) = \hbox{inv}_0(M_n(A))$,
with the result, the group
$\hbox{inv}(M_n(A))/$\break $\hbox{inv}_0(M_n(A))$ is isomorphic to the
group $\hbox{inv}(M_n(A_q))/\hbox{inv}_0(M_n(A_q))$. Hence by the
definition of the $K$-theory group $K_1$,
\begin{align*}
K_1(A) &= \lim\limits_\rightarrow
(\hbox{inv}(M_n(A))/\hbox{inv}_0(M_n(A)))\\[.2pc]
&= \lim\limits_\rightarrow
(\hbox{inv}(M_n(A_q))/\hbox{inv}_0(M_n(A_q))) = K_1(A_q).
\end{align*}
For $B$ to be $A$ or $A_q$, let the suspension of $B$ be
\begin{equation*}
SB = \{f\in C([0,1],B)\!: f(0) = f(1) = 0\} \cong C_0(\R, B).
\end{equation*}

We use the Bott periodicity theorem $K_0(B) = K_1(SB)$ to show
that $K_0(A) = K_0(A_q)$. It is standard that $\hbox{rad}(SA) =
\hbox{rad}(C_0(\R, A)) \cong C_0(\R, \hbox{rad}(A))$. Hence
\begin{align*}
SA_{q} &= C_{0} (\R, A_{q}) = C_{0} (\R, A/\hbox{rad}(A)) \cong
C_{0} (\R, A)/ C_{0} (\R, \hbox{rad}(A))\\[.2pc]
&= C_{0} (\R, A) /\hbox{rad} (C_{0} (\R, A)) = SA/\hbox{rad}(A)).
\end{align*}
Hence
\begin{equation*}
K_{0} (A_{q}) = K_{1}(SA_{q}) = K_{1} (SA/\hbox{rad}(SA)) =
K_{0}(A).
\end{equation*}
Thus we have
\begin{equation*}
K_{*}(A) = K_{*}(A_{q}) = K_{*}(E(A)) = RK_{*}(A) = RK_{*}(A_{q}).
\end{equation*}
Now $A^{\infty}$ is spectrally invariant in $A$ (Theorem~2.2 of
\cite{15}); and the action $\alpha$ on $A^{\infty}$ is smooth
(Theorem~A.2 of \cite{14}). Then applying the Phillips--Schweitzer
analogue of Thom isomorphism for smooth Frechet algebra crossed
product (Theorem~1.2 of \cite{11}) and Connes analogue of Thom
isomorphism for $C^{*}$-algebra crossed product \cite{7}, it
follows that
\begin{align*}
RK_{*} (S(\R, A^{\infty}, \alpha)) &= RK_{* + 1} (A^{\infty}) =
RK_{* + 1} (A) = RK_{* + 1} (E(A))\\[.2pc]
&=RK_{*} (C^{*} (\R, E(A), \alpha)) = K_{*}(C^{*}(\R, E(A),
\alpha)).
\end{align*}
When $\alpha$ is isometric, Theorem~1.3.4 of \cite{11} implies
that $RK_{*} (S(\R, A^{\infty}, \alpha)) = RK_{*}(L^{1}(\R, A,
\alpha))$. This completes the proof.\hfill$\Box$

\section{An application to the differential structure in
$\pmb{C^{*}}$-algebras}

Let $\cU$ be a unital $^{*}$-algebra. Let $\|\cdot\|$ be a
$C^{*}$-norm on $\cU$. Let $(A, \|\cdot\|)$ be the completion of
$(\cU, \|\cdot\|)$. Following \cite{5}, a map $T\hbox{\rm :}\ \cU
\rightarrow l^{1}(\N)$ is a {\it differential seminorm} if $T(x) =
(T_{k}(x))_{0}^{\infty} \in l^{1}(\N)$ satisfies the following:
\begin{enumerate}
\renewcommand\labelenumi{(\roman{enumi})}
\leftskip .4pc
\item $T_{k} (x) \geq 0$ for all $k$ and for all $x$.

\item For all $x, y$ in $\cU$ and scalars $\lambda, T(x + y) \leq
T(x) + T(y), T(\lambda x) = |\lambda|T(x)$.

\item For all $x, y$ in $\cU$, for all $k$,
\begin{equation*}
\hskip -1.25pc T_{k} (xy) \leq \sum_{i + j = k} T_{i} (x)
T_{j}(y).
\end{equation*}
\item There exists a constant $c > 0$ such that $T_{0}(x) \leq
c\|x\| \ \forall x \in \cU$.\vspace{.5pc}

By (ii), each $T_{k}$ is a seminorm. We say that $T$ is a {\it
differential} $^{*}$-{\it seminorm} if additionally;

\item $T_{k}(x^{*}) = T_{k}(x)$ for all $x$ and for all $k$.
\end{enumerate}

Further $T$ is a {\it differential norm} if $T(x) = 0$ implies $x
= 0$. Throughout we assume that $T_{0}(x) = \|x\|, x \in \cU$. The
{\it total norm} of $T$ is $T_{\rm tot} (x) = \sum_{k =
0}^{\infty} T_{k} (x), x \in \cU$. Given $T$, the differential
Frechet $^{*}$-algebra defined by $T$ is constructed as follows.
For each $k$, let $p_{k}(x) = \sum_{i = 0}^{k} T_{i} (x), x \in
\cU$. Then each $p_{k}$ is a submultiplicative $^{*}$-norm; and on
$\cU$, we have
\begin{equation*}
p_{0} \leq p_{1} \leq p_{2} \leq \cdots \leq p_{k} \leq p_{k + 1}
\leq \cdots
\end{equation*}
and $(p_{k})_{0}^{\infty}$ is a separating family of
submultiplicative $^{*}$-norms on $\cU$. Let $\tau$ be the locally
convex $^{*}$-algebra topology on $\cU$ defined by
$(p_{k})_{0}^{\infty}$. Let $\cU_{\tau} = (\cU, \tau)^{\sim}$ the
completion of $\cU$ in $\tau$ and let $\cU_{(k)} = (\cU,
p_{k})^{\sim}$ the completion of $\cU$ in $p_{k}$. Then
$\cU_{\tau}$ is a Frechet locally $m$-convex $^{*}$-algebra,
$\cU_{(k)}$ is a Banach $^{*}$-algebra. Let $\cU_{T}$ be the
completion of $(\cU, T_{\rm tot})$. Then the Banach $^{*}$-algebra
$\cU_{T} = \{x \in \cU_{\tau}\hbox{\rm :}\ \sup_{n} p_{n} (x) <
\infty\}$, the bounded part of $\cU_{\tau}$. By the definitions,
there exists continuous surjective $^{*}$-homomorphisms
$\phi_{k}\hbox{\rm :}\ \cU_{(k)} \rightarrow A, \phi\hbox{\rm :}\
\cU_{\tau} \rightarrow A$. The identity map $\cU \rightarrow \cU$
extends uniquely as continuous surjective $^{*}$-homomorphisms
$\varphi_{k}\hbox{\rm :}\ \cU_{(k + 1)} \rightarrow \cU_{(k)}$
such that
\begin{equation*}
{\cU}_{(0)} \xleftarrow{\ \ \ \ \varphi_{0}\ \ \ } {\cU}_{(1)}
\xleftarrow{\ \ \ \ \varphi_{1}\ \ \ } {\cU}_{(2)} \xleftarrow{\ \
\ \ \varphi_{2}\ \ \ } {\cU}_{(3)} \xleftarrow{\ \ \ \ \ \ \ \ \
}\cdots
\end{equation*}
is a dense inverse limit sequence of Banach $^{*}$-algebras and
$\cU_{\tau} = \lim\limits_{\leftarrow} \cU_{(k)}$.

\setcounter{theora}{0}
\begin{lem}\hskip -.3pc{\rm \cite{4}.}\ \
Let $(\cU,\|\cdot \|)$ be a $C^{*}$-normed algebra. Let $A$ be the
completion of \ $\cU$. Let $B$ denote $\cU_{(k)}$ or $\cU_{\tau}$
with respective topologies. Then the following hold{\rm :}
\begin{enumerate}
\renewcommand\labelenumi{{\rm (\roman{enumi})}}
\leftskip .4pc
\item $B$ is a hermitian $Q$-algebra.
\item $E(B) = A$.
\item $K_{*}(B) = K_{*}(A) = RK_{*} (B)$.
\end{enumerate}
\end{lem}
The $K$-theory result follows from the following.

\begin{lem}\hskip -.3pc{\rm \cite{4}.}\ \
Let $A$ be a Frechet algebra in which each element is bounded. Let
$A$ be spectrally invariant in $E(A)$. Then $K_{*} (A) = K_{*}
(E(A))$.
\end{lem}

Now let $\alpha$ be an action of $\R$ on $A$ leaving $\cU$
invariant. Let $T$ be $\alpha$-invariant, i.e., $T_{k}(\alpha(x))
= T_{k}(x)$ for all $k$ and for all $x$. Then $\alpha$ induces
isometric actions of $\R$ on each of $\cU_{(k)}, \cU_{\tau}$ and
$\cU_{T}$. Let $B$ be as above. Hence the crossed product Frechet
$^{*}$-algebras $L^{1}(\R, B^{\infty}, \alpha), L^{1}(\R, B,
\alpha), S(\R, B, \alpha)$ and $S(\R, B^{\infty}, \alpha)$ are
defined. Theorem~2 and Lemma~6.1 give the following, which is
Theorem~3(a).

\begin{coro}$\left.\right.$\vspace{.5pc}

\noindent $RK_{*}(S(\R, B^{\infty}, \alpha)) = RK_{*}(S(\R, B,
\alpha)) = RK_{*}(C^{*}(\R, A, \alpha)) = K_{* + 1} (A)$.
\end{coro}

Now let $\tilde{\cU}$ be the completion of $\cU$ in the family
$\cF$ of all $\alpha$-invariant differential $^{*}$-norms on
$\cU$. Then $\tilde{\cU}$ is a complete locally $m$-convex
$^{*}$-algebra admitting a continuous surjective
$^{*}$-homomorphism $\Psi \hbox{\rm :}\ \tilde{\cU} \rightarrow
A$. This $\alpha$-{\it invariant smooth envelope} $\tilde{\cU}$ is
different from the smooth envelope defined in \cite{5}, and it
need not be a subalgebra of $A$.

\begin{lem}
Assume that $\tilde{\cU}$ is metrizable. Then $\tilde{\cU}$ is a
hermitian Q-algebra{\rm ,} $E(\tilde{\cU}) = A${\rm ,} and $K_{*}
(\tilde{\cU}) = K_{*} (A)$.
\end{lem}

\looseness 1 This supplements a comment in p.~279 of \cite{5} that
$K_{*}(A) = \dot{K}_{*} (\cU_{1})$ where $\cU_{1}$ is the
completion of $\cU$ in all, not necessarily $\alpha$-invariant nor
closable, differential seminorms.

\begin{proof}
Since $\tilde{\cU} = \lim_{\leftarrow} \cU_{\tau}$, we have
$E(\tilde{\cU}) = \lim_{\leftarrow} E(\cU_{\tau}) = A$; and
$\tilde{\cU}$ admits greatest continuous $C^{*}$-seminorm, say
$p_{\infty} (\cdot)$ \cite{1}. It is easily seen that for any $x
\in \tilde{\cU}$, the spectral radius in $\tilde{\cU} r(x) \leq
p_{\infty} (x)$; and $\tilde{\cU}$ is a hermitian $Q$-algebra.
This implies, in view of $\tilde{E}(\tilde{\cU}) = A$, that the
spectrum in $\tilde{\cU} \hbox{sp}(x) = \hbox{sp}_{A}(j(x))$ for
all $x$ in $\tilde{\cU}$, where $j(x) = x + \hbox{srad}\
\tilde{\cU}$.

It follows from Lemma~6.2 that $K_{*}(A) = K_{*}(E(A))$. Hence
Lemma~6.4 follows.\hfill $\Box$
\end{proof}

Now the action $\alpha$ induces an isometric action of $\R$ on
$\tilde{\cU}$, with the result that the crossed product algebras
$S(\R, \tilde{\cU}, \alpha)$ and $L^{1}(\R, \tilde{\cU}, \alpha)$
are defined and are complete locally \hbox{$m$-convex} $^{*}$-algebras
with a $C^{*}$-enveloping algebras satisfying
\begin{align*}
E(S(\R, \tilde{\cU}, \alpha)) &= E(L^{1}(\R, \tilde{\cU},
\alpha))\\[.2pc]
&= C^{*} (\R, E(\tilde{\cU}), \alpha)\\[.2pc]
&= C^{*} (\R, A, \alpha).
\end{align*}
Theorem~2 quickly gives the following which is Theorem~3(b).

\begin{coro}$\left.\right.$\vspace{.5pc}

\noindent Assume that $\tilde{\cU}$ is metrizable. Then
$RK_{*}(S(\R, \tilde{\cU}, \alpha)) = K_{* + 1} (A)$.
\end{coro}

\end{document}